\documentclass{fpsac}

\usepackage[latin1]{inputenc}
\usepackage[T1]{fontenc}
\usepackage{ae,aecompl}

\usepackage{amssymb,amscd}
\usepackage{algorithmic}
\usepackage[matrix,arrow,curve]{xy}

\usepackage{verbatim}

\title[representation theories of some towers of algebras]
{Representation theories of some towers of algebras related to the
  symmetric groups and their Hecke algebras}
\date{15/11/2005}

\author{Florent Hivert and Nicolas M.~Thiéry}

\address{Florent Hivert -- \rm\texttt{florent.hivert@univ-rouen.fr}
 -- LIFAR -- university of Rouen --  Avenue de l'Université --
  Technopôle du Madrillet --
76801 Saint Etienne du Rouvray cedex -- FRANCE}

\address{Nicolas M. Thiéry -- \rm\texttt{nthiery@users.sf.net}
  -- Laboratoire de Mathématiques -- Université Paris Sud, Bât 425 --
  91405 Orsay Cedex -- FRANCE}



\newcommand{\TODO}[1]{}
\newcommand{\FIXME}[1]{}

\renewcommand{\k}{{\mathbb{C}}} % MAYBE CHANGE BACK TO \K?
\newcommand{\NN}{{\mathbb{N}}}
\newcommand{\F}{{\mathbb{F}}}
\newcommand{\CC}{{\mathbb{C}}}
\newcommand{\ZZ}{{\mathbb{Z}}}
\newcommand{\s}{\sigma}
\renewcommand{\t}{\tau}
\newcommand{\End}{{\operatorname{End}}}
\newcommand{\sym}{{\operatorname{Sym}}}
\newcommand{\qsym}{{\operatorname{QSym}}}
\newcommand{\ncsf}{{\operatorname{NCSF}}}

\newcommand{\sg}[1]{{\mathfrak{S}_{#1}}}
\newcommand{\sga}[2][\CC]{{#1[\sg{#2}]}}
\newcommand{\ksg}[2][\CC]{{#1\sg{#2}}}
\newcommand{\hecke}[1]{{\operatorname{H}_{#1}}}
\newcommand{\heckesg}[1]{{\operatorname{H\mathfrak{S}}_{#1}}}
\newcommand{\affinehecke}[1]{\operatorname{\hat H}_{#1}}
\newcommand{\ndf}[1]{\operatorname{NDF}_{#1}}
\newcommand{\ndfa}[2][\CC]{{#1[\ndf{#2}]}}
\newcommand{\ndpf}[1]{\operatorname{NDPF}_{#1}}
\newcommand{\ndpfa}[2][\CC]{{#1[\ndpf{#2}]}}

\newcommand{\Des}{{\operatorname{Des}}}
\newcommand{\Rec}{{\operatorname{Rec}}}
\newcommand{\cset}{{\operatorname{C}}}
\newcommand{\compof}{\vDash}                    % Composisition de
\newcommand{\fin}{{\succeq}}                    % Ordre de raffinement
\newcommand{\lon}{\operatorname{\ell}}          % Longeur
\newcommand{\act}{\cdot}

\newcommand{\idemp} {{p}}
\newcommand{\ls} {{\overline \s}}

\newcommand{\opi}{{\overline \pi}}

\newcommand{\init}{{\operatorname{init}}}
\newcommand{\id}{{\operatorname{id}}}

\newcommand{\sign}{{\operatorname{sign}}}

\newcommand{\rad}{{\operatorname{rad}}}
\newcommand{\cartan}{C}

\newcommand{\calG}{{\mathcal G}}
\newcommand{\calK}{{\mathcal K}}

\newtheorem{theorem}{Theorem}[section]
\newtheorem{lemma}{Lemma}[section]
\newtheorem{proposition}{Proposition}
\newtheorem{corollary}{Corollary}

\theoremstyle{definition}
\newtheorem{definition}{Definition}

\newtheorem{conjecture}{Conjecture}

\theoremstyle{remark}

\begin{document}

\begin{abstract}
  We study the representation theory of three towers of algebras which
  are related to the symmetric groups and their Hecke algebras. The
  first one is constructed as the algebras generated simultaneously by
  the elementary transpositions and the elementary sorting operators
  acting on permutations. The two others are the monoid algebras of
  nondecreasing functions and nondecreasing parking functions. For
  these three towers, we describe the structure of simple and
  indecomposable projective modules, together with the Cartan map. The
  Grothendieck algebras and coalgebras given respectively by the
  induction product and the restriction coproduct are also given
  explicitly.  This yields some new interpretations of the classical
  bases of quasi-symmetric and noncommutative symmetric functions as
  well as some new bases.

\begin{resume}
  Nous étudions la théorie des représentations de trois tours
  d'algèbres liées aux groupes symétriques et à leurs algèbres de
  Hecke. La première est formée des algèbres engendrées par les
  transpositions élémentaires ainsi que les opérateurs de tris
  élémentaires agissant sur les permutations. Les deux autres sont
  formées des algèbres des monoïdes des fonctions croissantes et des
  fonctions de parking croissantes. Pour ces trois tours, nous donnons
  la structure des modules simples et projectifs indécomposables ainsi
  que l'application de Cartan. Nous calculons également explicitement
  les algèbres et cogèbres de Grothendieck pour le produit d'induction
  et le coproduit de restriction. Il en découle de nouvelles
  interprétations de bases connues des fonctions quasi-symétriques et
  symétriques noncommutatives ainsi que des nouvelles bases.
\end{resume}
\end{abstract}

\maketitle

\tableofcontents

\section{Introduction}

Given an \emph{inductive tower of algebras}, that is a sequence of algebras
\begin{equation}
A_0 \hookrightarrow A_1 \hookrightarrow \cdots \hookrightarrow A_n
\hookrightarrow \cdots,
\end{equation}
with embeddings $A_m\otimes A_n \hookrightarrow A_{m+n}$ satisfying an
appropriate associativity condition, one can introduce two \emph{Grothendieck
rings}
\begin{equation}
\calG(A):=\bigoplus_{n\ge 0}G_0(A_n)\qquad \text{and} \qquad
\calK(A):=\bigoplus_{n\ge 0}K_0(A_n)\,,
\end{equation}
where $G_0(A)$ and $K_0(A)$ are the (complexified) Grothendieck groups
of the categories of finite-dimensional $A$-modules and projective
$A$-modules respectively, with the multiplication of the classes of an
$A_m$-module $M$ and an $A_n$-module $N$ defined by the induction product
\begin{equation}
[M] \cdot [N] = [M\widehat{\otimes} N] =
[M\otimes N \uparrow_{A_m\otimes A_n}^{A_{m+n}}]\,.
\end{equation}

If $A_{m+n}$ is a projective $A_m\otimes A_n$ modules, one can define a
coproduct on these rings by means of restriction of representations,
turning these into coalgebras. Under favorable circumstances the
product and the coproduct are compatible turning these into mutually
dual Hopf algebras.

The basic example of this situation is the character ring of the
symmetric groups (over $\CC$), due to Frobenius. Here the
$A_n:=\sga{n}$ are semi-simple algebras, so that
\begin{equation}
{G}_0(A_n) = {K}_0(A_n)= R(A_n)\,,
\end{equation}
where $R(A_n)$ denotes the vector space spanned by isomorphism classes
of indecomposable modules which, in this case, are all simple and
projective.
The irreducible representations $[\lambda]$ of $A_n$ are parametrized by
partitions $\lambda$ of $n$, and the Grothendieck ring is isomorphic to the
algebra $\sym$ of symmetric functions under the
correspondence
%\begin{equation}
$[\lambda] \leftrightarrow s_\lambda$,
%\end{equation}
where $s_\lambda$ denotes the Schur function associated with
$\lambda$.  Other known examples with towers of group algebras over
the complex numbers $A_n:=\CC[G_n]$ include the cases of wreath
products $G_n := \Gamma\wr\sg{n}$ (Specht), finite linear groups $G_n
:= GL(n,\F_q)$ (Green), \emph{etc.}, all related to symmetric
functions (see~\cite{Mcd,Zel}).

Examples involving non-semisimple specializations of Hecke algebras have also
been worked out.
Finite Hecke algebras of type $A$ at roots of unity ($A_n=H_n(\zeta)$,
$\zeta^r=1$) yield quotients and subalgebras of $Sym$~\cite{LLT}.
%% \begin{equation}
%% {\mathcal G} = Sym/(p_{rm}=0), \quad
%% {\mathcal K} = \CC\left[p_k\ |\ k\not\equiv0\pmod r\right]
%% \end{equation}
The Ariki-Koike algebras at roots of unity give rise to level $r$ Fock spaces
of affine Lie algebras of type $A$~\cite{AK}.
The $0$-Hecke algebras $A_n=\hecke{n}(0)$ correspond to the pair
Quasi-symmetric functions / Noncommutative symmetric functions,
${\mathcal G}=\qsym$, ${\mathcal K}=\ncsf$~\cite{NCSF4}.  Affine Hecke
algebras at roots of unity lead to $U(\widehat{sl}_r)$ and
$U(\widehat{sl}_r)^*$~\cite{Ari}, and the case of affine Hecke
generic algebras can be reduced to a subcategory admitting as
Grothendieck rings $U(\widehat{gl}_\infty)$ and
$U(\widehat{gl}_\infty)^*$~\cite{Ari}.
Further interesting examples are the tower of $0$-Hecke-Clifford
algebras~\cite{Ols,BHT} giving rise to the peak algebras~\cite{Stem},
and a degenerated version of the Ariki-Koike algebras~\cite{HNTAriki}
giving rise to a colored version of $\qsym$ and $\ncsf$.

\bigskip

The goal of this article is to study the representation theories of
several towers of algebras which are related to the symmetric groups
and their Hecke algebras $\hecke{n}(q)$. We describe their
representation theory and the Grothendieck algebras and coalgebras
arising from them. Here is the structure of the paper together with
the main results.

In Section~\ref{section.heckesg}, we introduce the main object of this
paper, namely a new tower of algebras denoted $\heckesg{n}$. Each
$\heckesg{n}$ is constructed as the algebra generated by both
elementary transpositions and elementary sorting operators acting on
permutations of $\{1,\dots,n\}$. We show that this algebra is better
understood as the algebra of antisymmetry preserving operators; this
allows us to compute its dimension and give an explicit basis.  Then,
we construct the projective and simple modules and compute their
restrictions and inductions. This gives rise to a new interpretation
of some bases of quasi-symmetric and noncommutative symmetric
functions in representation theory. The Cartan matrix suggests a link
between $\heckesg{n}$ and the incidence algebra of the boolean
lattice. We actually show that these algebra are Morita equivalent.
We conclude this section by discussing some links with a certain
central specialization of the affine Hecke algebra.

In Sections~\ref{section.ndf} and ~\ref{section.ndpf} we turn to the
study of two other towers, namely the towers of the monoids algebras
of nondecreasing functions and of nondecreasing parking functions.  In
both cases, we give the structure of projective and simple modules,
the cartan matrices, and the induction and restrictions rules.  We
also show that the algebra of nondecreasing parking functions is
isomorphic to the incidence algebra of some lattice.  Finally, we
prove that those two algebras are the respective quotients of
$\heckesg{n}$ and $\hecke n(0)$, through their representations on
exterior powers. The following diagram summarizes the relations
between all the mentioned towers of algebras:
\medskip
\begin{equation}
  \vcenter{
  \xymatrix@R=1cm@C=1cm{
    \hecke{n}(-1)        \ar@{->>}[d] \ar@{^(->}@/^4ex/[rrrr] &
    \hecke{n}(0)         \ar@{->>}[d] \ar@{^(->}@/^3ex/[rrr]  &
    \hecke{n}(1)=\sga{n} \ar@{->>}[d] \ar@{^(->}@/^2ex/[rr]   &
    \hecke{n}(q)         \ar@{->>}[d] \ar@{^(->}    [r]       &
    H\sg{n}              \ar@{->>}[d]\\
    \operatorname{Temperley-Lieb}_n                  \ar@{^(->}@/_4ex/[rrrr]&
    \ndpfa n                                         \ar@{^(->}@/_3ex/[rrr] &
    \sga{n}     \hookrightarrow \bigwedge^\cdot \k^n \ar@{^(->}@/_2ex/[rr]  &
    \hecke{n}(q)\hookrightarrow \bigwedge^\cdot \k^n \ar@{^(->}    [r]      &
    \ndfa n
  }}
\end{equation}

% l'action de \sg_n sur les produits exterieurs est une algebre de dimension
%        binom(2*n, n) = sum binom(n,i)^2

%\subsection*{Acknowledgment}

This paper mostly reports on a computation driven research using the
package \texttt{MuPAD-Combinat} by the authors of the present
paper~\cite{MuPAD-Combinat}. This package is designed for the computer
algebra system \texttt{MuPAD} and is freely available from
\texttt{http://mupad-combinat.sf.net/}. Among other things, it allows
to automatically compute the dimensions of simple and indecomposable
projective modules together with the Cartan invariants matrix of a
finite dimensional algebra, knowing its multiplication table.

\section{Background}

\subsection{Compositions and sets}

\TODO{Maybe take systematically $s_1,\dots,s_p$ as elements of $S$,
  and use $I=(i_1,\dots,i_p)$ here instead of $K$.}

Let $n$ be a fixed integer. Recall that each subset $S$ of $\{1,\dots,n-1\}$
can be uniquely identified with a $p$-tuple $K := (k_1,\dots,k_p)$ of
positive integers of sum $n$:
\begin{equation}
  \label{eq.set.to.comp}
  S=\{i_1 < i_2 < \dots < i_p\}
  \longmapsto \cset(S):=(i_1,i_2-i_1,i_3-i_2,\ldots,n-i_p)\,.
\end{equation}
We say that $K$ is a \emph{composition of $n$} and we write it by
$K\compof n$. The converse bijection, sending a composition to its
\emph{descent set}, is given by:
\begin{equation}
  \label{eq.comp.to.set}
  K = (k_1,\dots,k_p) \longmapsto
  \Des(K) = \{k_1+\dots+k_j,\  j=1,\dots,p-1\}\,.
\end{equation}
The number $p$ is called the \emph{length} of $K$ and is denoted by $\lon(K)$.

%% A weakly decreasing composition is called a \emph{partition}. A $p$-tuple of
%% non-negative integers of sum $n$ is called a \emph{pseudo composition}. Let
%% $I=(i_1,\ldots,i_q)$ and $J=(j_1,\ldots,j_p)$ be two compositions. By
%% $I\catcomp J$ we mean the concatenation of the two compositions $I\catcomp
%% J=(i_1,\ldots,i_q,j_1,\ldots,j_p)$. We denote by $I\sumcomp J$ the composition
%% defined by $(i_1,\ldots,i_q+j_1,\ldots,j_p)$. Sometimes we need the
%% composition of length $r$ with all parts equal to $i$; it will be denoted by
%% $(i^r)$ or briefly by $i^r$.
%% For instance, the composition $(3,1,2,1,2,2)$ of $11$ corresponds to
%% the subset $\{3,4,6,7,9\}$ of $\{1,2,3,4,5,6,7,8,9,10\}$.

The notions of complementary of a set $S^c$ and of inclusion of sets
can be transfered to compositions, leading to the complementary of a
composition $K^c$ and to the refinement order on compositions: we say
that $I$ is {\em finer} than $J$, and write $I\fin J$, if and only if
$\Des(I)\supseteq \Des(J)$.

%% For
%% example, the conjugate of $(3,2,1,4)$ is $(1,1,2,3,1,1,1)$. These two
%% compositions correspond to descent sets $\{3,5,6\}$ and $\{1,2,4,7,8,9\}$.
%% This can be read on compositions in the following way: Let
%% $J=(j_1,\ldots,j_p)$. The composition $I$ is finer than $J$ iff there exist
%% compositions $I_1$ of $j_1$ , $I_2$ of $j_2$, \ldots, $I_p$ of $j_p$ such that
%% $I=I_1\catcomp I_2\cdots I_p$ is the composition obtained by gluing
%% $I_1,\ldots,I_p$ one after another.

\subsection{Symmetric groups and Hecke algebras}

Take $n\in\NN$ and let $\sg{n}$ be the $n$-th symmetric group. It is
well known that it is generated by the $n-1$ elementary transpositions
$\sigma_i$ which exchange $i$ and $i+1$, with the relations
\begin{alignat}{2}
\s_i^2             & = 1                     && (1\leq i\leq n-1)\,, \notag\\
\s_i \s_j          & = \s_j \s_i             && (|i-j|\geq2)\,, \\
\s_i \s_{i+1} \s_i & = \s_{i+1} \s_i \s_{i+1} &\quad& (1\leq i\leq n-2)\,.\notag
\end{alignat}
The last two relations are called the \emph{braids relations}.  A
\emph{reduced word} for a permutation $\mu$ is a decomposition
$\mu=\sigma_{i_1}\cdots\sigma_{i_k}$ of minimal length. When denoting
permutations we also use the \emph{word notation}, where $\mu$ is
denoted by the word $\mu_1\mu_2\cdots\mu_n :=
\mu(1)\mu(2)\cdots\mu(n)$. For a permutation $\mu$, the set $\{i,\,
\mu_i > \mu_{i+1}\}$ of its \emph{descents} is denoted $\Des(\mu)$.
The descents of the inverse of $\mu$ are called the \emph{recoils of
  $\mu$} and their set is denoted $\Rec(\mu)$. For a composition $I$,
we denote by $\sg{I} := \sg{i_1} \times \dots \times \sg{i_p}$ the
\emph{standard Young subgroup} of $\sg{n}$, which is generated by the
elementary transpositions $\s_i$ where $i \notin \Des(I)$.

Recall that the (Iwahori-) \emph{Hecke algebra} $\hecke{n}(q)$ of type
$A_{n-1}$ is the $\CC$-algebra generated by elements $T_i$ for $i<n$
with the braids relations together with the quadratic relations:
\begin{equation}
  T_i^2=(q-1)T_i+q\,,
\end{equation}
where $q$ is a complex number.

The $0$-Hecke algebra is obtained by setting $q=0$ in these relations. Then,
the first relation becomes $T_i^2=-T_i$~\cite{Nor,NCSF4}. In this paper, we
prefer to use another set of generators $(\pi_i)_{i=1\dots n-1}$ defined by
$\pi_i := T_i + 1$. They also satisfy the braids relations together with the
quadratic relations $\pi_i^2 = \pi_i$.

Let $\sigma =: \sigma_{i_1}\cdots\sigma_{i_p}$ be a reduced word for a
permutation $\sigma\in\sg{n}$. The defining relations of $\hecke{n}(q)$
ensures that the element $T_{\sigma} := T_{i_1}\cdots T_{i_p}$ (resp.:
$\pi_{\sigma} := \pi_{i_1}\cdots \pi_{i_p}$) is independent of the chosen
reduced word for $\sigma$. Moreover, the well-defined family
$(T_{\sigma})_{\sigma\in\sg{n}}$ (resp.: $(\pi_{\sigma})_{\sigma\in\sg{n}}$)
is a basis of the Hecke algebra, which is consequently of dimension $n!$.

\subsection{Representation theory}

In this paper, we mostly consider \emph{right} modules over algebras.
Consequently the composition of two endomorphisms $f$ and $g$ is
denoted by $fg= g \circ f$ and their action on a vector $v$ is written
$v \act f$. Thus $g \circ f (v) = g(f(v))$ is denoted $v \act fg = (v
\act f) \act g$.

It is known that $\hecke{n}(0)$ has $2^{n-1}$ simple modules, all
one-dimensional, and naturally labelled by compositions $I$ of
$n$~\cite{Nor}: following the notation of~\cite{NCSF4}, let $\eta_I$
be the generator of the simple $\hecke{n}(0)$-module $S_I$ associated
with $I$ in the left regular representation. It satisfies
\begin{equation}
\label{etaI}
\eta_I \act T_i :=
\begin{cases}
  -\eta_I & \text{if $i\in\Des(I)$,}\\
  0       & \text{otherwise,}
\end{cases}
\qquad\text{or equivalently}\qquad
\eta_I \act \pi_i :=
\begin{cases}
  0       & \text{if $i\in\Des(I)$,}\\
  \eta_I  & \text{otherwise.}
\end{cases}
\end{equation}
The bases of the indecomposable projective modules $P_I$ associated to
the simple module $S_I$ of $\hecke{n}(0)$ are indexed by the
permutations $\sigma$ whose descents composition is $I$.

The Grothendieck rings of $\hecke{n}(0)$ are naturally isomorphic to
the dual pair of Hopf algebras of quasi-symmetric functions $\qsym$ of
Gessel~\cite{Ges} and of noncommutative symmetric functions
$\ncsf$~\cite{NCSF1} (see~\cite{NCSF4}). The reader who is not
familiar with those should refer to these papers, as we will only
recall the required notations here.

The Hopf algebra $\qsym$ of quasi-symmetric functions has two
remarkable bases, namely the \emph{monomial basis} $(M_I)_I$ and the
\emph{fundamental basis} (also called \emph{quasi-ribbon}) $(F_I)_I$.
They are related by
\begin{equation}
  F_I = \sum_{I\fin J} M_J\,
\qquad\text{or equivalently}\qquad
  M_I = \sum_{I\fin J} (-1)^{\lon(I) - \lon(J)} F_J\,.
\end{equation}
The characteristic map $S_I \mapsto F_I$ which sends the simple
$\hecke{n}(0)$ module $S_I$ to its corresponding fundamental function
$F_I$ also sends the induction product to the product of $\qsym$ and
the restriction coproduct to the coproduct of $\qsym$.

The Hopf algebra $\ncsf$ of noncommutative symmetric
functions~\cite{NCSF1} is a noncommutative analogue of the algebra of
symmetric functions~\cite{Mcd}. It has for multiplicative bases the
analogues $(\Lambda^I)_I$ of the elementary symmetric functions
$(e_\lambda)_\lambda$ and as well as the analogues $(S^I)_I$ of the
complete symmetric functions $(h_\lambda)_\lambda$. The relevant basis
in the representation theory of $\hecke{n}(0)$ is the basis of so
called \emph{ribbon Schur functions} $(R_I)_I$ which is an analogue of
skew Schur functions of ribbon shape. It is related to $(\Lambda_I)_I$
and $(S_I)_I$ by
\begin{equation}
  S_I = \sum_{I\fin J} R_J\,
  \qquad\text{and}\qquad
  \Lambda_I = \sum_{I\fin J} R_{J^c}\,.
\end{equation}
Their interpretation in representation theory goes as follows. The complete
function $S^n$ is the characteristic of the trivial module $S_{n} \approx
P_{n}$, the elementary function $\Lambda^n$ being the characteristic of the
sign module $S_{1^n} \approx P_{1^n}$. An arbitrary indecomposable projective
module $P_I$ has $R_I$ for characteristic. Once again the map $P_I \mapsto
R_I$ is an isomorphism of Hopf algebras.

Recall that $S_J$ is the semi-simple module associated to $P_I$,
giving rise to the duality between $\calG$ and $\calK$ :
\begin{equation}
  S_I = P_J/\rad(P_J)
  \qquad\text{and}\qquad
  \langle P_I\ ,\ S_J \rangle = \delta_{I,J}
\end{equation}
This translates into $\qsym$ and $\ncsf$ by setting that $(F_I)_I$ and
$(R_I)_I$ are dual bases, or equivalently that $(M_I)_I$ and $(S^I)_I$
are dual bases.

%%%%%%%%%%%%%%%%%%%%%%%%%%%%%%%%%%%%%%%%%%%%%%%%%%%%%%%%%%%%%%%%%%%%%%%%%%%%%%
\section{The algebra $\heckesg{n}$}
\label{section.heckesg}

The algebra of the symmetric group $\sga{n}$ and the $0$-Hecke algebra
$\hecke{n}(0)$ can be realized simultaneously as operator algebras by
identifying the underlying vector spaces of their right regular
representations.

Namely, consider the plain \emph{vector space} $\ksg{n}$
(distinguished from the \emph{group algebra} which is denoted by
$\sga{n}$).  On the first hand, the algebra $\sga{n}$ acts naturally
on $\ksg{n}$ by multiplication on the right (action on positions).
That is, a transposition $\sigma_i$ acts on a permutation $\mu :=
(\mu_1, \dots, \mu_n)$ by permuting $\mu_i$ and $\mu_{i+1}$:
$\mu\act\sigma_i=\mu\sigma_i$.

On the other hand, the $0$-Hecke algebra $\hecke{n}(0)$ acts on the
right on $\ksg{n}$ by decreasing sort. That is, $\pi_i$ acts on the
right on $\mu$ by:
\begin{equation}
  \mu \act \pi_i =
  \begin{cases}
    \mu          & \text{if $\mu_i > \mu_{i+1}$,}\\
    \mu \sigma_i & \text{otherwise}.
  \end{cases}
\end{equation}

\begin{definition}
  For each $n$, the algebra $\heckesg{n}$ is the subalgebra of
  $\End(\ksg{n})$ generated by both sets of operators
  $\s_1,\dots,\s_{n-1}, \pi_1,\dots,\pi_{n-1}$.
\end{definition}

By construction, the algebra $\heckesg{n}$ contains both $\sga{n}$ and
$\hecke{n}(0)$. In fact, it contains simultaneously all the Hecke
algebras $\hecke{n}(q)$ for all values of $q$; each one can be
realized by taking the subalgebra generated by the operators:
\begin{equation}
  T_i:=(q-1) (1-\pi_i) + q \s_i, \quad  \text{ for } i=1,\dots, n-1 \ .
\end{equation}

The natural embedding of $\ksg{n}\otimes\ksg{m}$ in $\ksg{n+m}$ makes
$(\heckesg{n})_{n\in\NN}$ into a tower of algebras, which contains the
similar towers of algebras $(\sga{n})_{n\in\NN}$ and
$(\hecke n(q))_{n\in\NN}$.

\subsection{Basic properties of $\heckesg{n}$}

Let $\opi_i$ be the \emph{increasing sort operator} on $\ksg{n}$.
Namely: $\opi_i$ acts on the right on $\mu$ by:
\begin{equation}
  \mu \act \opi_i =
  \begin{cases}
    \mu     & \text{ if $\mu_i < \mu_{i+1}$,}\\
    \mu \sigma_i & \text{otherwise.}
  \end{cases}
\end{equation}
Since $\pi_i + \opi_i$ is a symmetrizing operator, we have the identity:
\begin{equation}
  \pi_i + \opi_i = 1 + \s_i\ .
\end{equation}
It follows that the operator $\opi_i$ also belongs to $\heckesg{n}$.

The following identities are also easily checked:
\begin{equation}
  \label{equation.sigmapi}
\begin{aligned}
  \sigma_i\pi_i &= \pi_i\ ,   &  \sigma_i\opi_i &= \opi_i\ , \\
  \opi_i\pi_i   &= \pi_i\ ,   &  \pi_i\opi_i    &= \opi_i\ , \\
  \pi_i\sigma_i &= \opi_i\ ,  &  \opi_i\sigma_i &= \pi_i\ .
\end{aligned}
\end{equation}

A computer exploration suggests that the dimension of $\heckesg{n}$ is
given by the following sequence (sequence A000275 of the encyclopedia
of integer sequences~\cite{sloane}):
\begin{displaymath}
  1,1,3,19,211,3651,90921,3081513,136407699,7642177651,528579161353,44237263696473,
  \dots
\end{displaymath}
These are the numbers $h_n$ of pairs $(\s,\t)$ of permutations such
that $\Des(\s)\cap \Des(\t)=\emptyset$. Together with
Equation~(\ref{equation.sigmapi}), this leads to state the following
\begin{theorem}
  \label{theorem.heckesg.base}
  A vector space basis of $\heckesg{n}$ is given by the family of operators
  \begin{equation}
    B_n := \left\{ \s\pi_\t \ |\
        \Des(\s)\cap \Des(\t^{-1}) =\emptyset \right\}\ .
  \end{equation}
\end{theorem}

% First attempt by mean of abstract presentation

One approach to prove this theorem would be to find a presentation of
the algebra. The following relations are easily proved to hold in
$\heckesg{n}$:
%% \TODO{Revert all of these for right actions}
%% \begin{equation}
%% \begin{gathered} % Original version
%% \sigma_{i} \pi_{i+1} = \pi_{i}\pi_{i+1} +
%%    \pi_{i+1}\pi_{i} \sigma_{i+1}\sigma_{i}  - \pi_{i}\pi_{i+1}\pi_{i}\ , \\
%% \sigma_{i+1} \pi_{i} = \pi_{i+1}\pi_{i} + \pi_{i}\pi_{i+1}
%%    \sigma_{i}\sigma_{i+1} - \pi_{i}\pi_{i+1}\pi_{i}\ , \\
%% \sigma_1 \sigma_2 \pi_1 = \pi_2 \sigma_1 \sigma_2\ , \\
%% \sigma_1 \pi_2 \sigma_1 = \sigma_2 \pi_1 \sigma_2\ ,
%% \end{gathered}
%% \end{equation}
\TODO{The following relations follow from the last one below. Are they
  needed in the algorithm? $\pi_1 \sigma_2 \sigma_1 = \sigma_2\sigma_1
  \pi_2 \ $, $\sigma_1 \sigma_2 \pi_1 = \pi_2\sigma_1 \sigma_2$.}
\begin{equation}
  \label{equation.rewrite}
\begin{gathered}
\pi_{i+1} \sigma_{i} = \pi_{i+1}\pi_{i} +
   \sigma_{i}\sigma_{i+1} \pi_{i}\pi_{i+1} - \pi_{i}\pi_{i+1}\pi_{i}\ , \\
\pi_{i} \sigma_{i+1} = \pi_{i}\pi_{i+1} +
   \sigma_{i+1}\sigma_{i} \pi_{i+1}\pi_{i} - \pi_{i}\pi_{i+1}\pi_{i}\ , \\
\sigma_1 \pi_2 \sigma_1 = \sigma_2 \pi_1 \sigma_2\ ,
\end{gathered}
\end{equation}
and we conjecture that they generate all relations.
\begin{conjecture}
  A presentation of $\heckesg{n}$ is given by the defining relations
  of $\sga{n}$ and $\hecke{n}(0)$ together with the relations
  $\sigma_i\pi_i = \pi_i$ and of Equations~(\ref{equation.rewrite}).
\end{conjecture}
Using those relations as rewriting rules yields a straightening
algorithm which rewrites any expression in the $\s_i$'s and $\pi_i$'s
into a linear combination of the $\s\pi_\t$.  This algorithm seems, in
practice and with an appropriate strategy, to always terminate.
However we have no proof of this fact; moreover this algorithm is not
efficient, due to the explosion of the number and length of words in
intermediate results.

This is a standard phenomenon with such algebras. Their properties
often become clearer when considering their concrete representations
(typically as operator algebras) rather than their abstract
presentation. Here, theorem \ref{theorem.heckesg.base} as well as the
representation theory of $\heckesg{n}$ follow from its upcoming
structural characterization as the algebra of operators preserving
certain anti-symmetries.

% Second attempt by mean of concrete representation
\subsection{$\heckesg{n}$ as algebra of antisymmetry-preserving operators}

Let $\ls_i$ be the \emph{right operator} in $\End(\ksg{n})$ describing the
action of $s_i$ by multiplication \emph{on the left} (action on values),
namely $\ls_i$ is defined by
\begin{equation}
  \label{eq.def.lsi}
  \sigma \act \ls_i := \sigma_i\sigma\ .
\end{equation}
A vector $v$ in $\ksg{n}$ is \emph{left $i$-symmetric} (resp.
\emph{antisymmetric}) if $v \act \ls_i = v$ (resp. $v \act \ls_i = -v$). The
subspace of left $i$-symmetric (resp. antisymmetric) vectors can be
alternatively described as the image (resp. kernel) of the idempotent operator
$\frac12(1+\ls_i)$, or as the kernel (resp. image) of the idempotent operator
$\frac12(1-\ls_i)$.
\begin{theorem}
  \label{theorem.heckesg.characterization}%
  $\heckesg{n}$ is the subspace of $\End(\ksg{n})$ defined by the $n-1$
  \emph{idempotent sandwich equations}:
  \begin{equation}
    \frac12(1 - \ls_i) f \frac12(1 + \ls_i) = 0,
    \quad  \text{ for } i=1,\dots, n-1 \ .
  \end{equation}
  In other words, $\heckesg{n}$ is the subalgebra of those operators in
  $\End(\ksg{n})$ which preserve left anti-symmetries.
\end{theorem}
Note that, $\ls_i$ being self-adjoint, the adjoint algebra of
$\heckesg{n}$ satisfies the equations:
\begin{equation}
  \frac12(1 + \ls_i) f \frac12(1 - \ls_i) = 0;
\end{equation}
thus, it is the subalgebra of those operators in $\End(\ksg{n})$ which
preserve left symmetries. The symmetric group algebra has a similar
description as the subalgebra of those operators in $\End(\ksg{n})$
which preserve both left symmetries and antisymmetries.
\TODO{This is equivalent to commuting with its left action.
  $\hecke{n}(0)$ is also characterized by its commutation with its left
  action; this can also be seen as saying that it preserves
  "i-non-symmetries" (meaning all terms in the vector $v$ are sorted
  decreasingly at a given value).}

\begin{proof}
  The proof of theorem~\ref{theorem.heckesg.characterization} proceeds
  as follow. We first exhibit a triangularity property of the
  operators in $B_n$; this proves that they are linearly independent,
  so that $\dim \heckesg{n} \geq h_n$. Let $<$ be any linear extension
  of the right permutahedron order.
% or of the length weak order
Given an endomorphism $f$ of $\ksg{n}$, we order the rows and columns of
its matrix $M:=[f_{\mu\nu}]$ accordingly to $<$, and denote by
$\init(f) := \min\{ \mu, \exists \nu, f_{\mu\nu}\ne 0\}$ the index of
the first non zero row of $M$.
\begin{lemma}
  \label{lemma.independent}
  (a) Let $f:=\s \pi_\t$ in $B_n$. Then, $\init(f)=\t$, and
  \begin{equation}
    f_{\t\nu} =
    \begin{cases}
      1 & \text{if $\nu\in \sg{\Des(\t^{-1})} \s^{-1}$}\\
      0 & \text{otherwise}
    \end{cases}\ .
  \end{equation}

  (b) The family $B_n$ is free.
\end{lemma}
%\begin{proof}
%  (a) is a direct corollary of remark~\ref{remarks.factors} (b).
%
%  (b) follows by triangularity: the operator $\s \pi_\t$ has
%  coefficient $m_{\t\s^{-1}}=1$, whereas any other operator
%  $\s'\pi_\t'$ such that $D(\s')\cap R(\t) = \emptyset$, $\t'\leq \t$
%  and $\s'\ne\s$ has coefficient $m_{\t\s^{-1}}=0$.
%\end{proof}

Then, we note that $\heckesg{n}$ preserves all antisymmetries, because
its generators $\s_i$ and $\pi_i$ do. It follows that $\heckesg{n}$
satisfies the sandwich equations. We conclude by giving an explicit
description of the sandwich equations. Given an endomorphism $f$ of
$\ksg{n}$, denote by $(f_{\mu,\nu})_{\mu,\nu}$ the coefficients of its
matrix in the natural permutation basis. Given two permutations
$\mu,\nu$, and an integer $i$ in $\{1,\dots,n-1\}$, let
$R_{\mu,\nu,i}$ be the linear form:
\begin{equation}
  R_{\mu,\nu,i} :
  \begin{cases}
    \End(\ksg{n}) &\mapsto \k\\
    f            &\mapsto f_{\mu,\nu} + f_{s_i\mu,\nu} - f_{\mu,s_i\nu} + f_{s_i\mu,s_i\nu}
  \end{cases}
\end{equation}
Given a pair of permutations $\mu,\nu$ having at least one descent in
common, set $R_{\mu,\nu}=R_{\mu,\nu,i}$, where $i$ is the smallest
common descent of $\mu$ and $\nu$ (the choice of the common descent
$i$ is, in fact, irrelevant). Finally, let $R_n:=\{ R_{\mu,\nu},\,
\Des(\mu) \cap \Des(\nu) \ne \emptyset \}$.

\begin{lemma}
  \label{lemma.relations.independent}
  (a) If an operator $f$ in $\End{\ksg{n}}$ preserves
  $i$-antisymmetries, then $R_{\mu,\nu,i}(f)=0$ for any permutations
  $\mu$ and $\nu$.

  (b) The $n!^2-h_n$ linear relations in $R_n$ are linearly
  independent.
\end{lemma}
%% \begin{proof}
%%   (a) By linearity, we may assume without loss of generality that $f$
%%   is an operator obtained by composition of the generators $s_i,\pi_i$
%%   of $\heckesg$. By remark~\ref{remarks.factors} (d), $\mu.w=\nu$ or
%%   $\mu.w=s_i\nu$ if and only if $s_i\mu.w=\nu$ or $s_i\mu.w=s_i\nu$.
%%   Hence, $R_{\mu,\nu,i}(f)=0$.

%%   (b) Let $R_{\mu,\nu}=R_{\mu,\nu,i}$ be some relation in $R_n$, which
%%   we can be viewed as an $n!\times n!$ matrix. This matrix has a
%%   coefficient $1$ at position $\mu,\nu$. Since $i$ is both a descent
%%   of $\mu$ and $\nu$, $s_i\mu < \mu$, and $s_i\nu < \nu$; so the three
%%   other non zero coefficients are either strictly higher or more to
%%   the left of the matrix. Furthermore, no other relation in
%%   $R_{\mu,\nu}$ has a non-zero coefficient at position $\mu,\nu$.
%%   Hence, by triangularity the linear relations in $R_n$ are linearly
%%   independent.
%% \end{proof}

Theorems~\ref{theorem.heckesg.base}
and~\ref{theorem.heckesg.characterization} follow.
\end{proof}

\subsection{The representation theory of $\heckesg{n}$}

\subsubsection{Projective modules of $\heckesg{n}$}

Recall that $\heckesg{n}$ is the algebra of operators preserving left
antisymmetries. Thus, given $S\subset \{1,\dots, n-1\}$, it is natural
to introduce the $\heckesg{n}$-submodule $\bigcap_{i\in S} \ker
(1+\ls_i)$ of the vectors in $\ksg{n}$ which are $i$-antisymmetric for
all $i\in S$. For the ease of notations, it turns out to be better to
index this module by the composition associated to the
\emph{complementary set}; thus we define
\begin{equation}
  \label{eq.def.P.I}
  P_I := \bigcap_{i \notin \Des(I)} \ker (1+\ls_i)\, .
\end{equation}
The goal of this section is to prove that the family of modules
$(P_I)_{I\compof n}$ forms a complete set of representatives of the
indecomposable projective modules of $\heckesg{n}$.

The simplest element of $P_I$ is:
\begin{equation}
  v_I := \sum_{\nu \in \sg{I}} (-1)^{l(\nu)} \nu,
\end{equation}
One easily shows that
\begin{lemma}
  $v_I$ generates $P_I$ as an $\heckesg{n}$-module.
\end{lemma}

Given a permutation $\s$, let $v_\s:=v_{\Rec(\s)} \s$ (recall that $\Rec(\s) =
\Des(\s^{-1})$). Note that $\s$ is the permutation of minimal length appearing
in $v_\s$. By triangularity, it follows that the family $(v_\s)_{\s\in\sg{n}}$
forms a vector space basis of $\ksg{n}$. The usefulness of this basis comes
from the fact that
\begin{proposition}
  \label{proposition.base.PI}
  For any composition $I := (i_1,\dots,i_k)$ of sum $n$, the families
  \begin{equation}
    \{v_I \act \s\ |\ \s \in \sg{n},\ \Rec(\s) \, \cap\, \Des(I) = \emptyset\}
    \qquad \text{and} \qquad
    \{v_\s\        |\ \s \in \sg{n},\ \Rec(\s) \, \cap\, \Des(I) = \emptyset\}
  \end{equation}
  are both vector space bases of $P_I$; in particular, $P_I$ is of
  dimension $\frac{n!}{i_1!i_2!\dots i_k!}$.
\end{proposition}

Since $\sg{n}$ and $\hecke{n}(0)$ are both sub-algebras of $\heckesg{n}$,
the space $P_I$ is naturally a module over them. The following
proposition elucidates its structure.
\begin{proposition}
  \label{proposition.restriction.PI}
  Let $(-1)$ denote the sign representation of the symmetric group as
  well as the corresponding representation of the Hecke algebra
  $\hecke{n}(0)$ (sending $T_i$ to $-1$, or equivalently $\pi_i$ to
  $0$).
  \begin{itemize}
  \item[(a)] As a $\sg{n}$ module, $P_I \approx (-1)\uparrow_{\sg{I}}^{\sg{n}} $; its
    character is the symmetric function $e_I:=e_{i_1}\cdots e_{i_k}$.
  \item[(b)] As a $\hecke{n}(0)$ module, $P_I \approx
    (-1)\uparrow_{\hecke{I}(0)}^{\hecke{n}(0)}$; it is a projective
    module whose character is the noncommutative symmetric function
    $\Lambda^I:=\Lambda_{i_1}\cdots \Lambda_{i_k}$.
  \item[(c)] In particular the $P_I$'s are non isomorphic as $\hecke{n}(0)$-modules
    and thus as $\heckesg{n}$-modules.
  \end{itemize}
\end{proposition}
We are now in position to state the main theorem of this section.
\begin{theorem}
  \label{theorem.projective}
  For $\sigma\in \sg{n}$, let $\idemp_\s \in \End(\ksg{n})$ denote the projector
  on $\k v_\s$ parallel to $\oplus_{\t\ne \s} \k v_\t$. Then,
  \begin{itemize}
  \item[(a)] The ideal $\idemp_\s\heckesg{n}$ is isomorphic to $P_{\Rec(\s)} =
    P_{\Des(\s^{-1})} $ as an $\heckesg{n}$ module;
  \item[(b)] The idempotents $\idemp_\s$ all belong to $\heckesg{n}$;
    they give a maximal decomposition of the identity into orthogonal
    idempotents in $\heckesg{n}$;
  \item[(c)] The family of modules $(P_I)_{I\compof n}$ forms a
    complete set of representatives of the indecomposable projective
    modules of $\heckesg{n}$.
  \end{itemize}
\end{theorem}
\begin{proof}
  Item (a) is an easy consequence of Proposition
  \ref{proposition.base.PI}. To prove (b) one needs to check that
  $\idemp_\s$ belongs to $\heckesg{n}$. This is done by showing that it
  preserves left antisymmetries. Then, since the $\idemp_\s$'s give a
  maximal decomposition of the identity in $\End(\ksg{n})$, they are as
  well a maximal decomposition of the identity in $\heckesg{n}$.
  Finally, Item (c) follows from (a) and (b) and Item (c) of
  Proposition \ref{proposition.restriction.PI}.
\end{proof}

\subsubsection{Simple modules}

The simple modules are obtained as quotients of the projective modules
by their radical:
\begin{theorem}
  \label{theorem.heckesg.simple}
  The modules $S_I:=P_I / \sum_{J\subsetneq I} P_J$ form a complete
  set of representatives of the simple modules of $\heckesg{n}$.
  Moreover, the projection of the family $\{v_\s,\, \Rec(\s)=I\}$ in
  $S_I$ forms a vector space basis of $S_I$.
\end{theorem}
The modules $S_I$ are closely related to the projective modules of the
$0$-Hecke algebra:
\begin{proposition}
  The restriction of the simple module $S_I$ to $\hecke{n}(0)$ is an
  indecomposable projective module whose characteristic is the noncommutative
  symmetric function $R_{I^c}$.
\end{proposition}

\subsubsection{Cartan's invariants matrix and the boolean lattice}

We now turn to the description of the Cartan matrix. Let $p_I := p_\alpha$
where $\alpha$ is the shortest permutation such that $\Rec(\alpha) = I$
(this choice is in fact irrelevant).
\begin{proposition}
  \label{proposition.heckesg.cartan}
  Let $I$ and $J$ be two subsets of $\{1,\dots,n\}$. Then,
  \begin{equation}
    \dim Hom(P_I, P_J) = \dim p_I \heckesg{n} p_J=
    \begin{cases}
      1 & \text{if $I\subset J$,} \\
      0 & \text{otherwise.}
    \end{cases}
  \end{equation}
\end{proposition}
In other words, the Cartan matrix of $\heckesg{n}$ is the incidence
matrix of the boolean lattice. This suggests that there is a close
relation between $\heckesg{n}$ and the incidence algebra of the boolean
lattice. Recall that the \emph{incidence algebra} $\k[P]$ of a
partially ordered set $(P, \leq_P)$ is the algebra whose basis
elements are indexed by the couples $(u, v)\in P^2$ such that $u
\leq_P v$ with the multiplication rule
\begin{equation}
  \label{eq.incidence.algebra.mult}
  (u,v) \cdot (u',v') =
  \begin{cases}
    (u, v') & \text{if $v=u'$,} \\
    0       & \text{otherwise.}
  \end{cases}
\end{equation}
An algebra is called \emph{elementary} (or sometimes reduced) if its
simple modules are all one dimensional. Starting from an algebra $A$,
it is possible to get a canonical elementary algebra by the following
process. Start with a maximal decomposition of the identity $1=\sum_i
e_i$ into orthogonal idempotents.  Two idempotents $e_i$ and $e_j$ are
\emph{conjugate} if $e_i$ can be written as $a e_j b$ where $a$ and
$b$ belongs to $A$, or equivalently, if the projective modules $e_i A$
and $e_j A$ are isomorphic. Select an idempotent $e_c$ in each
conjugacy classes $c$ and put $e := \sum e_c$. Then, it is well
known~\cite{CR} that the algebra $eAe$ is elementary and that the
functor $M \mapsto M e$ which sends a right $A$ module to a $eAe$
module is an equivalence of category. Recall finally that two algebra
$A$ and $B$ such that the category of $A$-modules and $B$-modules are
equivalent are said \emph{Morita equivalent}.  Thus $A$ and $eAe$ are
Morita-equivalent.

Applying this to $\heckesg{n}$, one gets
\begin{theorem}
  Let $e$ be the idempotent defined by $e := \sum_{I \compof n} p_I$.
  Then the algebra $e \heckesg{n} e$ is isomorphic to the incidence
  algebra $\k[B_{n-1}]$ of the boolean lattice $B_{n-1}$ of subsets of
  $\{1,\dots,n-1\}$.  Consequently, $\heckesg{n}$ and $\k[B_{n-1}]$
  are Morita equivalent.
\end{theorem}

%% \begin{proof}
%%   \TODO{Rewrite this proof independently of the orthonormalization process}
%%   Let $f$ be an operator in $\End(\ksg{n})$, and $f'=p_I f p_J$. Then, the
%%   kernel of $f'$ contains $\ker p_I$ which is of codimension $1$,
%%   while the image of $f'$ is contained in $\im p_J$ which is of
%%   dimension $1$.  It follows that the sandwiches $p_I \End(\ksg{n}) p_J$,
%%   and $p_I \heckesg{n} p_J$ are of dimension at most $1$.

%%   Since $\im p_I=\k v_I$ and $v_I.\heckesg{n}=v_i.\k.\sg{n}$, we can
%%   furthermore conclude that $p_I \heckesg{n} p_J$ is of dimension $1$
%%   if and only if there exists a permutation $\s$ such that $v_I\s
%%   p_J\ne 0$. By proposition~\ref{proposition.idempotents}(b), this
%%   implies $I\subset J$.  Reciprocally, if $I\subset J$, taking
%%   $\s=\id$ yields $v_I\s p_J=v_J\ne 0$.
%% \end{proof}

\subsubsection{Induction, restriction, and Grothendieck rings}

Let $\calG := \calG\left((\heckesg{n})_n\right)$ and $\calK:=
\calK\left((\heckesg{n})_n\right)$ be respectively the Grothendieck
rings of the characters of the simple and projective modules of the
tower of algebras $(\heckesg{n})_n$.  Let furthermore $\cartan$ be the
cartan map from $\calK$ to $\calG$. It is the algebra and coalgebra
morphism which gives the projection of a module onto the direct sum of
its composition factors. It is given by
\begin{equation}
  \cartan(P_I) = \sum_{I\fin J} S_J\,.
\end{equation}
Since the indecomposable projective modules are indexed by
compositions, it comes out as no surprise that the structure of
algebras and coalgebras of $\calG$ and $\calK$ are each isomorphic to
$\qsym$ and $\ncsf$.  However, we do not get Hopf algebras, because
the structures of algebras and coalgebras are not compatible.

\begin{proposition}
  The following diagram gives a complete description of the structures
  of algebras and of coalgebras on $\calG$ and $\calK$.
  \begin{equation}
    \entrymodifiers={+<20pt>}%[F]}%{=<10pt>[o]}
    \vcenter{
    \xymatrix@R=0pt@C=1.8cm{
      (\qsym, .)
      & ({\calG},.)
      \ar@{^(->>}[l]_{\chi(S_I)\mapsto M_{I^c}}
      & ({\calK},.)
      \ar@{^(->>}[l]_{\cartan}
      \ar@{^(->>}[r]^{\chi(P_I)\mapsto F_{I^c}}
      & (\qsym, .)
      \\
      (\ncsf, \Delta)
      & ({\calG},\Delta) %\ar[ru]_{*}
      \ar@{^(->>}[l]^{\chi(S_I)\mapsto R_{I^c}}
      & ({\calK},\Delta) %\ar[lu]_{\not *}
      \ar@{^(->>}[l]^{\cartan}
      \ar@{^(->>}[r]_{\chi(P_I)\mapsto \Lambda^I}
      & (\ncsf, \Delta)
    }}
  \end{equation}
\end{proposition}
\begin{proof}
  The bottom line is already known from Proposition
  \ref{proposition.restriction.PI} and the fact that, for all $m$ and
  $n$, the following diagram commutes
\begin{equation}
  \entrymodifiers={+<10pt>}%[F]}%{=<10pt>[o]}
  \vcenter{\xymatrix{
      {\hecke{m}(0) \otimes \hecke{n}(0)\ } \ar@{^(->}[d] \ar@{^(->}[r] &
      {\ \hecke{m+n}(0)} \ar@{^(->}[d]
      \\
      {\heckesg{m} \otimes \heckesg{n}\ } \ar@{^(->}[r] &
      {\ \heckesg{m+n}}
    }}
\end{equation}
Thus the map which sends a module to the characteristic of its
restriction to $\hecke{n}(0)$ is a coalgebra morphism. The isomorphism
from $({\calK},.)$ to $QSym$ is then obtained by Frobenius duality
between induction of projective modules and restriction of simple
modules. And the last case is obtained by applying the Cartan map $\cartan$.
\end{proof}
It is important to note that the algebra $({\calG},.)$ is not the dual
of the coalgebra $(\calK, \Delta)$ because the dual of the restriction
of projective modules is the so called \emph{co-induction} of simple
modules which is, in general, not the same as the induction for non
self-injective algebras.

Finally the same process applied to the adjoint algebra which preserve
symmetries would have given the following diagram
  \begin{equation}
    \entrymodifiers={+<20pt>}%[F]}%{=<10pt>[o]}
    \vcenter{
    \xymatrix@R=0pt@C=1.8cm{
      (\qsym, .)
      & ({\calG},.)
      \ar@{^(->>}[l]_{\chi(S_I)\mapsto X_{I^c}}
      & ({\calK},.)
      \ar@{^(->>}[l]_{\cartan}
      \ar@{^(->>}[r]^{\chi(P_I)\mapsto F_{I}}
      & (\qsym, .)
      \\
      (\ncsf, \Delta)
      & ({\calG},\Delta) %\ar[ru]_{*}
      \ar@{^(->>}[l]^{\chi(S_I)\mapsto R_{I}}
      & ({\calK},\Delta) %\ar[lu]_{\not *}
      \ar@{^(->>}[l]^{\cartan}
      \ar@{^(->>}[r]_{\chi(P_I)\mapsto S^I}
      & (\ncsf, \Delta)
    }}
  \end{equation}
where $(X_I)_I$ is the dual basis of the elementary basis $(\Lambda_I)_I$ of
$\ncsf$. Thus we have a representation theoretical interpretation of many
bases of $\ncsf$ and $\qsym$.

%% The first step of the proof is to describe the restriction rule for
%% projective modules.
%% \begin{lemma}
%%   \TODO{coproduct of h or e in NCSF}
%% \end{lemma}
%% %% \begin{proof}
%% %%   $\sg{n+m}$ splits as $\sg{n}\times \sg{m}$ module into the direct sum
%% %%   $\bigoplus_{E\subset \{1,\dots,n+m\}} \sg{E} \otimes \sg{\overline
%% %%     E}$, which again split as expected.
%% %% \end{proof}
%% Then, all the other isomorphisms follow by applying the Frobenius
%% duality between restriction of projective modules and induction of
%% simple modules, and by using the $\cartan$ morphism which is
%% one-to-one (because the cartan matrix is upper triangular).

\subsection{Links with the affine Hecke algebra}
\newcommand{\affheckcent}[1]{{\mathcal H}_{#1}}

Recall that, for any complex number $q$, the extended affine Hecke algebra
$\affinehecke{n}(q)$ of type $A_{n-1}$ is the $\CC$-algebra generated by
$(T_i)_{i=1\cdots n-1}$ together with an extra generator $\Omega$ verifying
the defining relations of the Hecke algebra and the relation:
\begin{equation}
  \label{eq.def.affine.hecke}
  \Omega T_{i-1} = T_i\Omega
  \qquad \text{for $1\le i \leq n$}.
\end{equation}
The center of the affine Hecke algebra is isomorphic to the ring of
symmetric polynomials in some variables $\xi_1,\dots, \xi_n$ and it
can thus be specialized. Let us denote $\affheckcent{n}(q)$ the
specialization of the center $\affinehecke{n}(q)$ to the alphabet $1,
q, \dots q^{n-1}$. That is
\begin{equation}
  \affheckcent{n}(q) :=
  \affinehecke{n}(q) / \langle e_i(\xi_1,\dots,\xi_n) - e_i(1, q, \dots
q^{n-1})\ |\ i=1\dots n\rangle\,.
\end{equation}
It is well known that the simple modules $S_I$ of $\affheckcent{n}(q)$
are indexed by compositions $I$ and that their bases are indexed by
descent classes of permutations. Thus one expects a strong link
between $\heckesg{n}$ and $\affheckcent{n}(q)$. It comes out as
follows. Let $q$ be a generic complex number (i.e.: not 0 nor a root
of the unity). Sending $\Omega$ to $\s_1\s_2\cdots\s_{n-1}$ and $T_i$
to itself yields a surjective morphism from $\affheckcent{n}(q)$ to
$\heckesg{n}$. Thus, the simple modules of $\affheckcent{n}(q)$ are
the simple modules of $\heckesg{n}$ lifted back through this morphism.
This also explains the link between the projective modules of
$\hecke{n}(0)$ and the simple modules of $\affheckcent{n}(q)$, thanks
to Proposition~\ref{proposition.restriction.PI}.

%%%%%%%%%%%%%%%%%%%%%%%%%%%%%%%%%%%%%%%%%%%%%%%%%%%%%%%%%%%%%%%%%%%%%%%%%%%%%%
\section{The algebra of non-decreasing functions}
\label{section.ndf}

\TODO{Maybe we should put a complement on the representations of
  $\ndfa n$ on sets, so that the fact that the projective modules start
  with k=1 coincides with the fact that the representations of
  $\heckesg n$ are indexed by subsets of size $<n$.}

\begin{definition}
  Let $\ndf{n}$ be the set of \emph{non-decreasing functions} from
  $\{1,\dots,n\}$ to itself. The composition and the neutral element
  $\id_n$ make $\ndf{n}$ into a monoid. Its cardinal is
  $\binom{2n-1}{n-1}$, and we denote by $\ndfa{n}$ its monoid algebra.
\end{definition}
The monoid $\ndf{n}\times \ndf{m}$ can be identified as the submonoid of
$\ndf{n+m}$ whose elements stabilize both $\{1,\dots,n\}$ and
$\{n+1,\dots,n+m\}$. This makes $(\ndfa{n})_n$ into a tower of algebras.

%% \FIXME{Maybe move this to a later place}
%% Given a function $f\in \ndf{n}$, we denote by $\im f$ the image of $f$,
%% $\rank f$ the size of $\im f$, and by $\ker f$ the set partition of
%% $\{1,\dots,n\}$ into $f$ fibers. Note that those fibers are
%% consecutive intervals, so that $\ker f$ is completely described by the
%% composition $I=(i_1,\dots,i_{\rank f})$, where $i_j$ is the size of
%% the $j$-th fiber. Furthermore, $f$ is completely determined by the
%% knowledge of $I$ and $\im f$.

One can take as generators for $\ndf{n}$ and $A_n$ the functions
$\pi_i$ et $\opi_i$, such that $\pi_i(i+1)=i$, $\pi_i(j)=j$ for $j\ne
i+1$, $\opi_i(i)=i+1$, and $\pi_i(j)=j$ for $j\ne i$. The functions
$\pi_i$ are idempotents, and satisfy the braid relations, together
with a new relation:
\begin{equation}
  \pi_i^2=\pi_i
  \qquad \text{and} \qquad
  \pi_{i+1}\pi_i\pi_{i+1} = \pi_i\pi_{i+1}\pi_i = \pi_{i+1}\pi_i\, .
\end{equation}
This readily defines a morphism $\phi: \pi_{\hecke{n}(0)} \mapsto
\pi_{\ndfa n}$ of $\hecke{n}(0)$ into $\ndfa n$. Its image is the
monoid algebra of \emph{non-decreasing parking functions} which will
be discussed in Section~\ref{section.ndpf}. The same properties hold
for the operators $\opi_i$'s.  Although this is not a priori obvious,
it will turn out that the two morphisms
$\phi:\pi_{\hecke{n}(0)}\mapsto \pi_{\ndfa n}$ and $\overline \phi:
\opi_{\hecke{n}(0)}\mapsto\opi_{\ndfa n}$ are compatible, making
$\ndfa n$ into a quotient of $\heckesg{n}$.

\subsection{Representation on exterior powers}

We now want to construct a suitable representation of $\ndfa{n}$ where
the existence of the epimorphism from $\heckesg{n}$ onto $\ndfa{n}$, and
the representation theory of $\ndfa{n}$ become clear.

The \emph{natural representation} of $\ndfa{n}$ is obtained by taking the
vector space $\k^n$ with canonical basis $e_1,\dots,e_n$, and letting a
function $f$ act on it by $e_i.f=e_{f(i)}$. For $n>2$, this representation is
a faithful representation of the monoid $\ndf{n}$ but not of the algebra, as
$\dim \ndfa{n}=\binom{2n-1}{n-1} \gg n^2$. However, since $\ndf{n}$ is a
monoid, the diagonal action on \emph{exterior powers}
\begin{equation}
  \label{eq.action.exterieur}
  (x_1\wedge\dots\wedge x_k) \act f :=
  (x_1 \act f)\wedge\dots\wedge (x_k \act f)
\end{equation}
%
%% \TODO{shall we say that $\ndfa{n}$ is a Hopf algebra? I don't think it is
%% useful to add the hopf formalism to this simple fact. If you agree with my
%% proposition, just remove this comment}%
%
still define an action. Taking the \emph{exterior powers} $\bigwedge^k
\k^n$ of the natural representation gives a new representation, whose
basis $\{e_S:=e_{s_1}\wedge\dots\wedge e_{s_k}\}$ is indexed by
subsets $S:=\{s_1,\dots,s_k\}$ of $\{1,\dots,n\}$. The action of a
function $f$ in $\ndf{n}$ is simply given by (note the absence of
sign!):
\begin{equation}
  e_S.f :=
  \begin{cases}
    e_{f(S)} & \text{ if $|f(S)|=|S|$},\\
    0      & \text{ otherwise.}
  \end{cases}
\end{equation}
We call \emph{representation of $\ndfa{n}$ on exterior powers} the
representation of $\ndfa{n}$ on $\bigoplus_{k=1}^n \bigwedge^k \k^n$,
which is of dimension $2^n-1$ (it turns out that we do not need to
include the component $\bigwedge^0 \k^n$ for our purposes).

\begin{lemma}
  \label{lemma.wifa.faithfull}
  The representation of $\ndfa{n}$ on $\bigoplus_{k=1}^n \bigwedge^k
  \k^n\bigwedge\k^n$ is faithful.
\end{lemma}

We now want to realize the representation of $\ndfa{n}$ on the $k$-th
exterior power as a representation of $\heckesg{n}$. To this end, we
use a variation on the standard construction of the Specht module
$V_{k,1,\dots,1}$ of $\sg{n}$ to make it a $\heckesg{n}$-module. The
trick is to use an appropriate quotient of $\k\sg{n}$ to simulate the
symmetries that we usually get by working with polynomials, while
preserving the $\heckesg{n}$-module structure. Namely, consider the
following $\heckesg{n}$-module:
\begin{equation}
  P_n^k := P_{k,1,\dots,1} / \bigcup P_{k,1,\dots,1,2,1,\dots,1}.
\end{equation}
An element in $P_n^k$ is left-antisymmetric on the values
$1,\dots,k-1$ and symmetric on the values $k+1,\dots,n-1$, the effect
of the quotient being to identify two permutations which differ by a
permutation of the values $\{k+1,\dots,n\}$. A basis of $P_n^k$
indexed by subsets of size $k$ of $\{1,\dots,n\}$ is obtained by
taking for each such subset $S$ the image in the quotient $P_n^k$ of
\begin{equation}
  e_S := \sum_{\sigma, \sigma(S) = \{1,\dots,k\},
    \sigma(i) < \sigma(j) \text{ for $i<j\not\in S$}}
  (-1)^{\sign \sigma} \sigma\, .
\end{equation}

\TODO{Generalize this construction to any Specht module?}

It is straightforward to check that the actions of $\pi_i$ and
$\opi_i$ of $\heckesg{n}$ on $e_S$ of $P_k$ coincide with the actions
of $\pi_i$ and $\opi_i$ of $\ndfa n$ on $e_S$ of $\bigwedge^k \k^n$
(justifying a posteriori the identical notations). In the sequel, we
identify the modules $P_n^k$ and $\bigwedge^k \k^n$ of $\heckesg{n}$
and $\ndfa n$, and we call \emph{representation on exterior powers of
  $\heckesg{n}$} its representation on $\bigoplus_{k=1}^n \bigwedge^k
\k^n$. Using Lemma~\ref{lemma.wifa.faithfull} we are in position to
state the following
\begin{proposition}
  $\ndfa{n}$ is the quotient of $\heckesg{n}$ obtained by considering
  its representation on exterior powers. The restriction of this
  representation of $\heckesg{n}$ to $\sga{n}$, $\hecke{n}(0)$, and
  $\hecke{n}(-1)$ yield respectively the usual representation of
  $\sg{n}$ on exterior powers, the algebra of non-decreasing parking
  functions (see Section~\ref{section.ndpf}), and the Temperley-Lieb
  algebra.
\end{proposition}

\TODO{blah blah blah, graded version of the natural representation on sets?}

\subsection{Representation theory}
\label{section.reptheo.ndf}

\subsubsection{Projective modules, simple modules, and Cartan's invariant matrix}
Let $\delta$ be the usual homology border map:
\begin{equation}
  \delta:
  \begin{cases}
    P_n^k & \to P_n^{k-1} \\
    S:=\{s_1,\dots,s_k\} & \mapsto \sum_{i\in \{1,\dots,k\}}
    (-1)^{k-i} S \backslash \{s_i\}
  \end{cases}\, .
\end{equation}
This map is naturally a morphism of $\ndfa n$-module. For each $k$ in
$1,\dots, n$, let $S_k := P_k / \ker \delta$. It turns out that
together with the identity, $\delta$ is essentially the only $\ndfa
n$-morphism. We are now in position to describe the projective and
simple modules, as well as the Cartan matrix of $\ndfa n$.
\begin{proposition}
  The modules $(P_n^k)_{k=1,\dots,n}$ form a complete set of
  representatives of the indecomposable projective modules of $\ndfa
  n$.

  The modules $(S_n^k)_{k=1,\dots,n}$ form a complete set of
  representatives of the simple modules of $\ndfa n$.

  \label{proposition.wifa.cartan}
  Let $k$ and $l$ be two integers in $\{1,\dots,n\}$. Then,
  \begin{equation}
    \dim Hom(P_n^k, P_n^l) =
    \begin{cases}
      1 & \text{if $l\in \{k, k-1\}$,} \\
      0 & \text{otherwise.}
    \end{cases}
  \end{equation}
\end{proposition}

The proof relies essentially on the following lemma:
\begin{lemma}
  There exists a minimal decomposition of the identity of $\ndfa n$
  into $2^n-1$ orthogonal idempotents. In particular, the
  representation on exterior powers is the smallest faithful
  representation of $\ndfa n$.
\end{lemma}

\TODO{Description of left and right modules in $\ndfa n$.}

\subsubsection{Induction, restriction, and Grothendieck groups}

\begin{proposition}
  The restriction and induction of indecomposable projective modules
  and simple modules are described by:
  \begin{equation}
    P_{n_1+n_2}^k \downarrow^{\ndfa{n_1+n_2}}_{\ndfa{n_1}\otimes \ndfa{n_2}} \approx
    \bigoplus_{\substack{
        n_1+n_2=n\\
        k_1+k_2=k\\
        1\leq k_i\leq n_i \text{ or } k_i=n_i=0
      }}
    P_{n_1}^{k_1} \otimes P_{n_2}^{k_2}
  \end{equation}
  \begin{equation}
    P_{n_1}^{k_1} \otimes P_{n_2}^{k_2} \uparrow^{\ndfa{n_1+n_2}}_{\ndfa{n_1}\otimes \ndfa{n_2}}\approx
    P_{n_1+n_2}^{k_1+k_2} \oplus P_{n_1+n_2}^{k_1+k_2-1}
  \end{equation}
  \begin{equation}
    S_{n_1+n_2}^k \downarrow^{\ndfa{n_1+n_2}}_{\ndfa{n_1}\otimes \ndfa{n_2}} =
    \bigoplus_{
      \substack{
        n_1+n_2=n\\
        k_1+k_2 \in \{k,k+1\}\\
        1\leq k_i\leq n_i \text{ or } k_i=n_i=0
      }
    }
    S_{n_1}^{k_1} \otimes S_{n_2}^{k_2}
  \end{equation}
  \begin{equation}
    S_{n_1}^{k_1} \otimes S_{n_2}^{k_2} \uparrow^{\ndfa{n_1+n_2}}_{\ndfa{n_1}\otimes \ndfa{n_2}}\approx
    S_{n_1+n_2}^{k_1+k_2}
  \end{equation}
\end{proposition}

Those rules yield structures of commutative algebras and cocommutative
coalgebras on $\calG$ and $\calK$ which can be realized as quotients
or sub(co)algebras of $\sym$, $\qsym$, and $\ncsf$.  However, we do
not get Hopf algebras, because the structures of algebras and
coalgebras are not compatible (compute for example
$\Delta(\chi(P_1^1)\chi(P_1^1))$ in the two ways, and check that the
coefficients of $\chi(P_1^1)\otimes\chi(P_1^1)$ differ).
%% \begin{proposition}
%%   The following diagram gives a complete description of the structure
%%   of algebra and coalgebra of $\calG$ and $\calK$.
%%   \begin{displaymath}
%%     \entrymodifiers={+<20pt>}%[F]}%{=<10pt>[o]}
%%     \xymatrix@R=0pt@C=1.8cm{
%%       (\sym, .)
%%       \ar@{->>}[r]^{h_\lambda \mapsto \chi\left(S_{|\lambda|}^{\lon(\lambda)}\right)}
%%       & (\calG,.)
%%       & (\calK,.)
%%       \ar@{^(->>}[l]_{\cartan}
%%       & (\ncsf, .)
%%       \ar@{->>}[l]^{R_I \mapsto \chi\left(P_{|I|}^{\lon(I)}\right)}
%%       \\
%%       (\qsym, \Delta)
%%       \ar@{->>}[r]^{F_I \mapsto \chi\left(S_{|I|}^{\lon(I)}\right)}
%%       & (\calG,\Delta) %\ar[ru]_{*}
%%       & (\calK,\Delta) %\ar[lu]_{\not *}
%%       \ar@{^(->>}[l]^{\cartan}
%%       \ar@{^(->}[r]_{\chi(P_n^k)\mapsto \sum_{\lambda\partof n,\lon(\lambda)=k} m_\lambda}
%%       & (\sym, \Delta)
%%     }
%%   \end{displaymath}
%% \end{proposition}

%%%%%%%%%%%%%%%%%%%%%%%%%%%%%%%%%%%%%%%%%%%%%%%%%%%%%%%%%%%%%%%%%%%%%%%%%%%%%%
\section{The algebra of non-decreasing parking functions}
\label{section.ndpf}

\begin{definition}
  A \emph{nondecreasing parking function} of size $n$ is a
  nondecreasing function $f$ from $\{1,2,\dots n\}$ to $\{1,2,\dots
  n\}$ such that $f(i) \leq i$, for all $i\leq n$.
  
  The composition of maps and the neutral element $\id_n$ make the set
  of nondecreasing parking function of size $n$ into a monoid denoted
  $\ndpf{n}$.
\end{definition}

It is well known that the nondecreasing parking functions are counted
by the Catalan numbers $C_n = \frac1{n+1}\binom{2n}{n}$. It is also
clear that $\ndpf{n}$ is the sub-monoid of $\ndf{n}$ generated by the
$\pi_i$'s.

\subsection{Simple modules}

The goal of the sequel is to study the representation theory of
$\ndpf{n}$, or equivalently of its algebra $\ndpfa{n}$. The following
remark allows us to deduce the representations of $\ndpfa{n}$ from the
representations of $\hecke{n}(0)$.
\begin{proposition}
  The kernel of the algebra epi-morphism $\phi : \hecke{n}(0) \to \ndpfa{n}$
  defined by $\phi(\pi_i) = \pi_i$ is a sub-ideal of the radical of
  $\hecke{n}(0)$.
\end{proposition}
\begin{proof}
  It is well known (see \cite{Nor}) that the quotient of $\hecke{n}(0)$ by its
  radical is a commutative algebra. Consequently, $\pi_i\pi_{i+1}\pi_i -
  \pi_i\pi_{i+1} = [\pi_i\pi_{i+1}, \pi_i]$ belongs to the radical of
  $\hecke{n}(0)$.
\end{proof}
As a consequence, taking the quotient by their respective radical shows that
the projection $\phi$ is an isomorphism from $\ndpfa{n}/\rad(\ndpfa{n})$ to
$\hecke{n}(0)/\rad(\hecke{n}(0))$. Moreover $\ndpfa{n}/\rad(\ndpfa{n})$ is
isomorphic to the commutative algebra generated by the $\pi_i$ such that
$\pi_i^2=\pi_i$. As a consequence, $\hecke{n}(0)$ and $\heckesg{n}$ share,
roughly speaking, the same simple modules:

\begin{corollary}
  There are $2^{n-1}$ simple $\ndpfa{n}$-modules $S_I$, and they are
  all one dimensional. The structure of the module $S_I$, generated by
  $\eta_I$, is given by
\begin{equation}
\left\{
\begin{array}{r@{\,}cl}
\eta_I\act\pi_i &= 0      & \text{if $i\in\Des(I)$,}\\
\eta_I\act\pi_i &= \eta_I & \text{otherwise.}
\end{array}
\right.
\end{equation}
\end{corollary}

\subsection{Projective modules}
The projective modules of $\ndpf{n}$ can be deduced from the ones of
$\ndf{n}$.
\begin{theorem}
  Let $I$ be a composition of $n$, and $S := \Des(I) = \{s_1, \dots,
  s_k\}$ be its associated set. Then, the principal \FIXME{translate :
    monogene} sub-module
  \begin{equation}
    P_I :=
    (e_1 \wedge e_{s_1+1} \wedge \dots \wedge e_{s_1+1}) \act \ndpfa{n}
    \quad \subset\quad \bigwedge^{k+1} \k^n
  \end{equation}
  is an indecomposable projective module. Moreover, the set $(P_I)_{I
    \compof n}$ is a complete set of representatives of indecomposable
  projective modules of $\ndpfa{n}$.
\end{theorem}

This suggests an alternative description of the algebra $\ndpfa{n}$.
Let $G_{n,k}$ be the lattice of subsets of $\{1,\dots,n\}$ of size $k$
for the \emph{product order} defined as follows. Let $S :=
\{s_1<s_2<\dots<s_k\}$ and $T := \{t_1<t_2<\dots<t_k\}$ be two
subsets.  Then,
\begin{equation}
  S \leq_G T
  \qquad\text{if and only if}\qquad
  s_i \leq t_i \text{, for $i=1,\dots, k$.}
\end{equation}
One easily sees that $S \leq_G T$ if and only if there exists a
nondecreasing parking function $f$ such that $e_S = e_T \act f$. This
lattice appears as the Bruhat order associated to the Grassman
manifold $G^n_k$ of $k$-dimensional subspaces in $\CC^n$.
\begin{theorem}
  There is a natural algebra isomorphism
  \begin{equation}
    \ndpfa{n}\ \approx\ \bigoplus_{k=0}^{n-1} \k[G_{n-1, k}]\ .
  \end{equation}
\end{theorem}
In particular the Cartan map $\cartan : \calK \to \calG$ is given by
the lattice $\leq_G$:
\begin{equation}
  \cartan(P_I)\ =\ \sum_{J,\ \Des(J) \leq_G \Des(I)} S_J
\end{equation}
On the other hand, due to the commutative diagram
\begin{equation}
  \entrymodifiers={+<10pt>}%[F]}%{=<10pt>[o]}
  \vcenter{\xymatrix{
      {\hecke{m}(0) \otimes \hecke{n}(0)\ } \ar@{->>}[d] \ar@{^(->}[r] &
      {\ \hecke{m+n}(0)} \ar@{->>}[d]
      \\
      {\ndpf{m} \otimes \ndpf{n}\ } \ar@{^(->}[r] &
      {\ \ndpf{m+n}}
    }}
\end{equation}
it is clear that the restriction of simple modules and the induction
of indecomposable projective modules follow the same rule as for
$\hecke{n}(0)$. The induction of simple modules can be deduced via the
Cartan map, giving rise to a new basis $G_I$ of $\ncsf$. It remains
finally to compute the restrictions of indecomposable projective
modules. It can be obtained by a not yet completely explicit
algorithm. All of this is summarized by the following diagram:
  \begin{equation}
    \entrymodifiers={+<20pt>}%[F]}%{=<10pt>[o]}
    \vcenter{
    \xymatrix@R=0pt@C=1.8cm{
      (\ncsf, .)
      & ({\calG},.)
      \ar@{^(->>}[l]_{\chi(S_I)\mapsto G_{I}}
      & ({\calK},.)
      \ar@{^(->>}[l]_{\cartan}
      \ar@{^(->>}[r]^{\chi(P_I)\mapsto R_{I}}
      & (\ncsf, .)
      \\
      (\qsym, \Delta)
      & (\calG,\Delta)
      \ar@{^(->>}[l]^{\chi(S_I)\mapsto F_I}
      & (\calK,\Delta)
      \ar@{^(->>}[l]^{\cartan}
      \ar@{^(->>}[r]_{\chi(P_I)\mapsto ???}
      & ???
    }}
\end{equation}

\end{document}